\documentclass{amsart}
\usepackage{graphicx,amssymb,stmaryrd,amsfonts,epsfig,amsthm,a4,amsmath,mathrsfs,url}
\usepackage{vmargin,amscd}
\usepackage{eufrak}
\usepackage[latin1]{inputenc}
\usepackage[all]{xy}
\psfigdriver{dvips}

\newtheorem{thm}{Theorem}
\newtheorem{cor}[thm]{Corollary}
\newtheorem{lem}[thm]{Lemma}

\newtheorem{prop?}[thm]{Proposition?}
\theoremstyle{definition}

\theoremstyle{remark}
\newtheorem*{remm}{Remark}
\newtheorem{rem}[thm]{Remark}

\numberwithin{equation}{section}

\newcommand{\Z}{\mathbf{Z}}

\newcommand{\R}{\mathbf{R}}
\newcommand{\C}{\mathbf{C}}

\newcommand{\tpr}{\begin{tiny}\noindent Proof:}

\newcommand{\bpr}{\noindent \textbf{Proof}: ~}

\newcommand{\epr}{~$\blacksquare$}

\begin{document}
\title{Semisimple Zariski closure of Coxeter groups}
\author{Yves de Cornulier}%
\date{July 19, 2005 (Minor revisions: Nov. 23, 2012)}





\begin{abstract}

Let $W$ be an irreducible, finitely generated Coxeter group. The
geometric representation provides an discrete embedding in the
orthogonal group of the so-called Tits form. One can look at the
representation modulo the kernel of this form; we give a new proof
of the following result of Vinberg: if $W$ is non-affine, then
this representation remains faithful. Our proof uses relative
Kazhdan Property (T).

The following
corollary was only known to hold when the Tits form is
non-degenerate: the reduced $C^*$-algebra of $W$ is simple with a
unique normalized trace. Some other corollaries are pointed
out.\end{abstract}

\maketitle





\setcounter{section}{-1}
\section{Foreword}
I wrote the present paper in the beginning of 2005. After Luis Paris informed me that its only purportedly original contribution was an old result of Vinberg, I left the paper on my web page as an unpublished note. I now post it on arXiv so as to make this access perennial.


\section{Introduction}
We recall that a (discrete) group $\Gamma$ is amenable if it has a
left invariant finitely additive probability on its power set. All
we need to know about the class of amenable group is that it is
stable under taking subgroups, quotient, and direct limits. It
immediately follows that every group $\Gamma$ has a unique biggest
amenable normal subgroup, which we call its \textit{amenable
radical} and denote it~$R(\Gamma)$.

Let $S$ be a set, and $M=(m_{st})_{(s,t)\in S^2}$ a $S\times S$
symmetric matrix, with $1$'s on the diagonal and coefficients in
$\{2,3,\dots,\infty\}$ outside the diagonal. The group $W$ with
presentation $$\big\langle (\sigma_s)_{s\in
S}\mid\;((\sigma_s\sigma_t)^{m_{st}})_{(s,t)\in S^2}\big\rangle,$$ where
(we set $x^\infty=1$ for all $x$) is called the Coxeter group
associated to the Coxeter matrix~$M$.

The Coxeter matrix $M$ defines a non-oriented labelled graph
$\mathcal{G}$ with $S$ as set of vertices and an edge between $s$
and $t$ if and only if $m_{s,t}\ge 3$; this edge being labelled by
$m_{st}$. The decomposition of the graph $\mathcal{G}$ in
connected components corresponds to a decomposition of $W$ into a
direct sum. Thus, most problems about the group structure of $W$,
such as the determination of $R(W)$, are reduced to the case when
$M$ is irreducible, which means by definition that $\mathcal{G}$
is connected and non-empty. By abuse of language, we sometimes say that $W$ is
irreducible.

Irreducible Coxeter groups fall into three classes:

\begin{itemize}
\item Coxeter groups of spherical type. They are locally finite;
they are entirely classified: the corresponding diagrams are
called $A_n$, ($n\ge 1$) $B_n$ ($n\ge 2$), $D_n$ ($n\ge 3$), $E_n$
$(6\le n\le 8)$, $F_4$, $H_3$, $H_4$, $I_2(n)$ ($n\ge 3$) for the
finitely generated (hence finite) ones, and $A_\infty$,
$A'_\infty$, $B_\infty$, $D_\infty$ for the infinitely generated
ones.

\item Coxeter groups of affine type. They are finitely generated,
infinite, virtually abelian. The corresponding diagrams are called
$\widetilde{A_n}$, ($n\ge 1$) $\widetilde{B_n}$ ($n\ge 3$),
$\widetilde{C_n}$ ($n\ge 2$) $\widetilde{D_n}$ ($n\ge 4$),
$\widetilde{E_n}$ $(6\le n\le 8)$, $\widetilde{F_4}$,
$\widetilde{G_2}$.

\item Non-affine Coxeter groups: these are the remaining Coxeter
groups. They contains a non-abelian free group \cite{Ha}. In the finitely generated case, they are
even \textit{large}, i.e.\ have a finite index subgroup mapping
onto a free non-abelian free group \cite{Gon, MaVi}. This is the
class we mainly study.
\end{itemize}


The purpose of the present paper is to exhibit a discrete Zariski
dense embedding of every non-affine Coxeter group in a simple
orthogonal group $\text{O}(p,q)$ for some $p\ge 2$, $q\ge 1$. This
is actually a combination of a result of Benoist and de la Harpe,
and a result of Vinberg, which we provide a new proof based on
Kazhdan's relative Property~(T).
Before this, let us begin by the corollaries.

\begin{cor}
Let $W$ be the Coxeter group associated to a non-affine
irreducible Coxeter matrix. Then its amenable radical is trivial: 
$R(W)=\{1\}$.\label{cor:radical_cox}
\end{cor}

We immediately deduce the amenable radical of an arbitrary Coxeter group.

\begin{cor}
Let $W$ be a Coxeter group. Let $W=\bigoplus W_i\oplus\bigoplus
W_j$ be its decomposition as a direct sum of irreducible Coxeter
groups, with $W_i$ non-affine and $W_j$ affine or spherical. Then
$R(W)=\bigoplus W_j$.
\end{cor}


Using a Theorem of Y. Benoist and P. de la Harpe \cite{BH}, which
also plays an essential role in this paper, 
we obtain several other consequences.

A group $\Gamma$ is called \textit{primitive} if it has a proper
maximal subgroup $\Lambda$ containing no nontrivial normal
subgroup of~$\Gamma$, i.e.\ such that the action of $\Gamma$ on
$\Gamma/\Lambda$ is faithful. Using a characterization of
primitive linear groups due to T.~Gelander and Y.~Glasner, we
obtain:


\begin{cor}
Let $W$ be an infinite, finitely generated Coxeter group. Then $W$
is primitive if and only if it is irreducible and non-affine.
\label{cor:cox_pro_dense}
\end{cor}


We also obtain a characterization of Coxeter groups with simple
reduced $C^*$-algebra. If $\Gamma$ is any discrete group and
$g\in\Gamma$, its regular representation of $\Gamma$ on
$\ell^2(\Gamma)$ is defined by $\lambda(g)(f)(h)=f(g^{-1}h)$ for
$f\in\ell^2(\Gamma)$. The norm closure of the linear span of the
operators $\lambda(g)$ for $g\in\Gamma$ is called the reduced
$C^*$-algebra of $\Gamma$ and denoted $C^*_r\Gamma$.


\begin{cor}
Let $W$ be a Coxeter group. Then its reduced $C^*$-algebra
$C^*_rW$ is simple if and only if $W$ has only non-affine factors.
Moreover, if these conditions are satisfied, then it has a unique
normalized trace.\label{cor:C_star_simple}
\end{cor}

The existence of a nontrivial amenable normal subgroup is a
well-known obstruction for $C^*_rW$ to be simple, and is
conjectured to be the only one for linear groups \cite{BekH}.
Actually Corollary \ref{cor:C_star_simple} reduces to show that
any non-affine irreducible Coxeter group has simple reduced
$C^*$-algebra. Actually this had been proved in \cite{Fe} and
\cite[§2, vii]{BH} under the additional assumption that the Tits
form (introduced below) is non-degenerate; in which case Corollary
\ref{cor:radical_cox} is an immediate consequence of results
of~\cite{BH}.

\begin{remm}
An extensive discussion on simplicity of group $C^*$-algebras can be found in \cite{Ha2}, where a part of Corollary \ref{cor:C_star_simple} (in the finitely generated irreducible case) appears as \cite[Cor.~18]{Ha2}.
\end{remm}

\medskip

Let us now recall the definition of the Tits form. On the vector
space $\R^{(S)}$ with basis $(e_s)_{s\in S}$, define a symmetric
bilinear form, called the \textit{Tits form}, by
$B(e_s,e_t)=-\cos(\pi/m_{st})$. Suppose now that~$S$ is finite,
and let $(p,q,r)$ be the signature of $B$, meaning that $B$ is
equivalent to the form~$\begin{pmatrix}
  I_p &  &  \\
   & -I_q &  \\
   &  & 0I_r \\
\end{pmatrix}$.

By abuse of language, we call $(p,q,r)$ the signature of~$W$. It
is known \cite[Chap. V, §4.8]{Bou} that the Coxeter system is
spherical if and only if $q=r=0$, i.e.\ $B$ is a scalar product. If
the Coxeter system is irreducible, then $q=0$ implies $r\le 1$,
and the Coxeter system is affine exactly when $(p,q,r)=(p,0,1)$
\cite[Chap. V, §4.9]{Bou}.

For $s\in S$, set $r_s(v)=v-2B(e_s,v)e_s$. The mapping
$\sigma_s\to r_s$ extends to a well-defined group morphism
$\alpha$, called the {\em Tits representation} of the Coxeter group
$W$. By a well-known theorem of Tits, this representation is
faithful and has discrete image \cite[Chap. V, §4.4]{Bou}.

Denote by $\text{O}_f(B)$ the subgroup of $\textnormal{GL}\big(\R^{(S)}\big)$
consisting of those linear maps preserving the form $B$ \textit{and}
fixing pointwise the kernel $\text{Ker}(B)$. It is known \cite[Chap.
V, §4.7]{Bou} that $\alpha(W)$ is contained in $\text{O}_f(B)$.

\begin{thm}[Benoist, de la Harpe \cite{BH}]
Suppose that $S$ is finite, $W$ is irreducible and non-affine. Then the image $\alpha(W)$ of the Tits representation is
Zariski-dense in $\text{O}_f(B)$.\label{thm:BH}
\end{thm}

The group $\text{O}_f(B)$ is easily seen to be isomorphic to
$\text{O}(p,q)\ltimes(\R^{p+q})^{r}$. Accordingly, when $r=0$, it
is isomorphic to $\text{O}(p,q)$, whose amenable radical is its
centre $\{\pm 1\}$. Since infinite irreducible Coxeter groups have
trivial centre \cite[Chap. V, §4, Exercice 3]{Bou}, one thus obtains
that if $r=0$, then the Coxeter group $W$ has a trivial amenable
radical.

However, even if, in a certain sense, ``most" Coxeter groups have
a non-degenerate Tits form, those with degenerate Tits form are
numerous and a classification of those seems out of reach. Here
are some examples:


For all $a,b,c,d\in\{2,3,4,\dots,\infty\}$ with $c,d\ge 3$, the
signature of the Coxeter diagram


$$
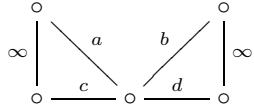
\xymatrix {\circ\ar^a@{-}[dr]\ar_\infty@{-}[d] &  &
\circ\ar^\infty@{-}[d]\\
\circ\ar^c@{-}[r] & \circ\ar^d@{-}[r]\ar^b@{-}[ur] & \circ}$$

\noindent is $(3,1,1)$. The signature of the Coxeter diagram


$$
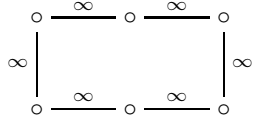
\xymatrix {\circ\ar^\infty@{-}[r]\ar_\infty@{-}[d] & \circ\ar^\infty@{-}[r] &
\circ\ar^\infty@{-}[d]\\
\circ\ar^\infty@{-}[r] & \circ\ar^\infty@{-}[r] & \circ}$$

\noindent is $(3,1,2)$. In \cite[§5]{BH}, irreducible Coxeter
groups with signature $\Big(p,1,(p+1)(p-2)/2\Big)$ are exhibited for all
$p\ge 4$.

\medskip

The bilinear form $B$ naturally factors through a non-degenerate
bilinear map $\bar{B}$ on the quotient $\R^{(S)}/\text{Ker}(B)$. Denote by $T_f$ the kernel of the natural map from
$\text{O}_f(B)$ to $\textnormal{O}(\bar{B})$, so that $T_f$ is
isomorphic to $(\R^{p+q})^r$. 
We have the following lemma (recall
that the Coxeter system is supposed finite, non-affine, and
irreducible):

\begin{lem}[Vinberg {\cite[Proposition 13]{V}}; see also {\cite[Proposition 6.1.3]{Kra}}]
$$\alpha(W)\cap T_f=\{1\}.$$\label{lemma_not_intersect_rad}
\end{lem}

A sketch of our proof of Lemma \ref{lemma_not_intersect_rad} goes as follows: set
$N=\alpha^{-1}(\alpha(W)\cap T_f)$. We must prove that $N=\{1\}$.
We prove that $(W,N)$ has relative Property (T), i.e.\ that every
isometric action of $W$ on a Hilbert space has a $N$-fixed point.
On the other hand, all finitely generated Coxeter groups are known
\cite{BoJS} to act properly by isometries on some Hilbert space
(this is called the Haagerup Property); the combination of these
two facts implies that $N$ is finite, hence trivial.


Lemma \ref{lemma_not_intersect_rad} has the following consequence. Let $\pi$ be the composite natural projection $\text{O}_f(B)\to\textnormal{O}(\bar{B})\to\textnormal{PO}(\bar{B})$.

\begin{thm}[Reduced Tits representation]
If $W$ is a non-affine irreducible finitely generated Coxeter group, then the
homomorphism $\pi\circ\alpha$ embeds $W$ as a discrete, Zariski
dense subgroup of~$\textnormal{PO}(\bar{B})$.\label{thm:discrete_proj}
\end{thm}

\begin{cor}\label{cn}
If $W$ is a non-affine irreducible finitely generated Coxeter group, then $W$ embeds as a discrete, $\C$-Zariski
dense subgroup of a complex simple Lie group with trivial centre, namely a projective orthogonal group.
\end{cor}

\medskip

In Section \ref{sec:all_but_lemma}, we show how Lemma
\ref{lemma_not_intersect_rad} implies all results above, and in
Section \ref{sec:proof_lemma} we prove Lemma
\ref{lemma_not_intersect_rad}. The reader may find superfluous to
use Kazhdan's Property (T) since a geometric proof already exists;
but we have included it here to illustrate a surprising
application of relative Property~(T).

\medskip

\textbf{Remerciements.} Je remercie Yves Benoist et Pierre de la
Harpe pour les discussions à ce sujet, et particulièrement Luis
Paris pour m'avoir signalé le résultat de Vinberg. Je remercie Goulnara Arzhantseva pour l'intérêt porté à ce papier.




\section{Proof of all results from Lemma \ref{lemma_not_intersect_rad}}\label{sec:all_but_lemma}

We use the easy lemma in the proof of Corollary
\ref{cor:radical_cox}:

\begin{lem}
If $S$ has at least 3 elements, and the signature of the
corresponding Coxeter group is $(p,q,r)$, then $p\ge
2$.\label{lem:pge2}
\end{lem}
\bpr If $S$ has an edge $st$ with finite label $m_{st}$ (possibly
$m_{st}=2$), then the restriction of $B$ to the plane $\R
e_s\oplus\R e_t$ is positive definite, so that $p\ge 2$. Otherwise
$S$ is the complete graph with all labels infinite, and, when $S$
has 3 elements, a direct computation shows that the signature is
$(2,1,0)$.\epr

\begin{rem}
By a less trivial result by Luis Paris \cite{Paris}, if $S$ has at
least 4 elements and is connected, then $p\ge 3$. On the other
hand we have pointed out above that an hexagon with infinite
labels has signature $(3,1,2)$; we do not know if there exist
irreducible Coxeter groups with $|S|\ge 7$ and
$p=3$.\label{rem:Paris}
\end{rem}

\medskip

\noindent \textbf{Proof of Theorem \ref{thm:discrete_proj}.} 
The injectivity of the map into $\textnormal{O}(\bar{B})$ is the contents of Lemma \ref{lemma_not_intersect_rad}. Since infinite irreducible Coxeter groups have trivial centre \cite[Chap. V, §4, Exercice 3]{Bou}, the composite map into $\textnormal{PO}(\bar{B})$ is still injective.

The Zariski density of its image follows from Theorem \ref{thm:BH}. It remains
to prove that the image $\Gamma=p\circ\alpha(W)$ is discrete. By a
theorem of Auslander \cite[Theorem 8.24]{Ragh}, the connected
component $\overline{\Gamma}^0$ (in the ordinary topology) is
solvable. Since $\Gamma$ is Zariski dense in the simple group
$\text{PO}(\bar{B})$, it follows that its normal subgroup
$\overline{\Gamma}^0$ must be discrete, hence trivial, i.e.\ $\Gamma$ is
discrete.\epr

\medskip

\noindent \textbf{Proof of Corollary \ref{cn}.}
By Theorem \ref{thm:discrete_proj}, $W$ embeds as a Zariski dense subgroup in the
real group $\text{PO}(\bar{B})\simeq \text{PO}(p,q)$ with $p\ge 2$ and $q\ge 1$. Taking the
complexification, we obtain, unless $(p,q)=(2,2)$ or $(3,1)$, a $\C$-Zariski
dense embedding of $W$ in $\text{PO}_{p+q}(\C)$ and since $4\neq p+q\ge 3$, we are done. If $(p,q)=(3,1)$, the group $\text{PO}(p,q)$ itself has a structure of a simple complex Lie group so the argument works without complexification. Finally, we cannot have $(p,q)=(2,2)$. Indeed, by a result of Paris (see Remark \ref{rem:Paris}), if $|S|\ge 4$ then $p\ge 3$.\epr

\medskip

\noindent \textbf{Proof of Corollary \ref{cor:radical_cox}.} We can suppose $W$ irreducible. If $W$ is finitely generated, then the result immediately follows from Corollary \ref{cn}. If $W$ is infinite, then the Coxeter graph is a direct limit of finite connected non-affine Coxeter subgraphs, and thus $W$ is a direct limit of finitely generated non-affine irreducible Coxeter groups. Since the property of having trivial amenable radical is stable under passing to direct limits, we are done.\epr


\medskip

\noindent \textbf{Proof of Corollary \ref{cor:cox_pro_dense}.}
It is a easy fact that if a primitive group $\Gamma$ decomposes as
a nontrivial direct product $\Gamma_1\times \Gamma_2$, then
$\Gamma_1$ and $\Gamma_2$ are simple non-abelian and isomorphic.
It immediately follows that any primitive Coxeter group $W$ must
be irreducible.

If $W$ is an affine Coxeter group, then it cannot be primitive:
indeed, let $M$ be a maximal subgroup. Then $W$ is virtually
abelian, hence is subgroup separable, i.e.\ every finitely
generated subgroup is the intersection of subgroups of finite
index containing it. Moreover every subgroup is finitely
generated. It immediately follows that every maximal subgroup in
$W$ must have finite index. Therefore, the action of the infinite
group $W$ on the finite set $W/M$ cannot be faithful.

Let us now suppose that $W$ is irreducible and non-affine.
Gelander and Glasner \cite{GG} prove that an infinite finitely generated
linear group $\Gamma$ is primitive if and only if there exists an
algebraically closed field $K$, a linear algebraic $K$-group $G$,
and a morphism $\Gamma\to G(K)$ with Zariski dense image, such
that $G^0$ is semisimple with trivial centre, and the action of
$\Gamma$ by conjugation on $G^0(K)$ is faithful and is transitive
on simple factors of $G^0$. It follows from Corollary \ref{cn} that this criterion is satisfied.\epr


\begin{rem}
The primitive finite Coxeter groups are those of type $A_n$ ($n\ge
1$), $D_{2n+1}$ ($n\ge 1$), $E_6$, and $I_2(p)$ for $p$ odd prime
(note the redundancies $I_2(3)\simeq A_2$, $A_3\simeq D_3$).
Indeed, Coxeter groups of type $B_n$ ($n\ge 2$), $D_{2n}$ ($n\ge
2$), $E_7$, $E_8$, $F_4$, $H_3$, $H_4$, $I_2(2n)$ ($n\ge 2$) have
centre cyclic of order 2. The only remaining cases are those of
type $I_2(n)$ for odd non-prime $n$, for which the verification is
straightforward.

Conversely, the group of type $A_n$, the symmetric group
$S_{n+1}$, acts primitively and faithfully on $n+1$ elements. The
group of type $D_{2n+1}$, isomorphic to
$S_{n+1}\ltimes(\Z/2\Z)^{n}_0$ acts affinely and primitively on
\[(\Z/2\Z)^{n}_0=\big\{(a_1,\dots,a_{n+1})\in(\Z/2\Z)^{n+1},
a_1+\dots+a_{n+1}=0\big\}.\] Finally, the group $W$ of type $E_6$ has
only one normal subgroup other than trivial ones, namely $W^+$, and it follows
that if $\Lambda$ is any maximal subgroup other than $W^+$, then
the action of $W$ on $W/\Lambda$ is faithful.
\end{rem}


\medskip

\noindent \textbf{Proof of Corollary \ref{cor:C_star_simple}.} If
$C^*_rW$ is simple, then $W$ has no nontrivial amenable normal
subgroup, so that $W$ has no affine or spherical factor.

Conversely, suppose that $W$ has no non-affine factor. Since the
property of having simple reduced $C^*$-algebra with a unique
normalized trace is inherited by direct limits \cite[Lemma
5.1]{BCH}, we are reduced to the case when $W$ is finitely
generated.

The argument is the same as that given in \cite[§2, vii]{BH},
except that, thanks to Lemma \ref{lemma_not_intersect_rad}, we can
avoid assuming that $B$ is non-degenerate.

We use the following results:

\begin{itemize}\item \cite[Theorem 1]{BCH} If a discrete group $\Gamma$ embeds as a Zariski dense
subgroup in a connected real semisimple Lie group without compact
factors, then $C^*_r\Gamma$ is simple and has a unique normalized
trace. \item \cite{BekH} If $\Gamma_0$ has finite index in
$\Gamma$, if $\Gamma$ is i.c.c.\ (all its nontrivial conjugacy
classes are infinite), and if $C^*_r\Gamma_0$ is simple and has a
unique normalized trace, then $C^*_r\Gamma$ is also simple and has
a unique normalized trace.
\end{itemize}

First suppose that $W$ is irreducible. By Corollary \ref{cor:radical_cox}, $W$ is i.c.c.; moreover $W$ has a subgroup of
index at most 2 embedding as a Zariski dense subgroup in
$\text{PO}_0(p,q)$. So the two criteria above apply.


In general, decompose $W$ as $W=W_1\times\dots\times W_n$ with each
$W_i$ irreducible non-affine. Then $W$ is i.c.c., and has a subgroup of index $\le 2^n$ that embeds as a
Zariski dense subgroup in a connected semisimple Lie group with
$n$ noncompact simple factors. It follows that $C^*_rW$ is also
simple and has a unique normalized trace.\epr

\begin{rem}
It is not difficult to see that, conversely, Lemma
\ref{lemma_not_intersect_rad} follows from any one among
Corollaries \ref{cor:radical_cox}, \ref{cor:cox_pro_dense}, or
~\ref{cor:C_star_simple}.
\end{rem}



\section{Proof of Lemma \ref{lemma_not_intersect_rad}}\label{sec:proof_lemma}

By \cite{BoJS}, if $W$ is any Coxeter group and $l$ its length
function, then there exists an isometric action $u$ of $W$ on a
Hilbert space $\mathcal{H}$, and $v\in\mathcal{H}$ such that
$l(g)=\|u(g)v-v\|^2$ for all $g\in W$.

On the other hand, recall that, given a group $\Gamma$ and a
subgroup $\Lambda$, the pair $(\Gamma,\Lambda)$ has
\textit{relative Property (T)} if for every isometric action $u$
of $\Gamma$ on a Hilbert space $\mathcal{H}$, and every
$v\in\mathcal{H}$, the restriction to $\Lambda$ of the function
$g\mapsto\|u(g)v-v\|$ is bounded.

It follows that if $W$ is a finitely generated Coxeter group, and
if $\Lambda\subset W$ is a subgroup such that $(W,\Lambda)$ has
relative Property (T), then $\Lambda$ is finite. In particular, if
$\Lambda$ is torsion-free, this implies $\Lambda=\{1\}$. Thus
Lemma \ref{lemma_not_intersect_rad} follows from the following
lemma:


\begin{lem}
Set $N=\alpha^{-1}(\alpha(W)\cap T_f)$. The pair $(W,N)$ has
relative Property (T).\label{lemma_rad_Trel}
\end{lem}

\bpr It is clear from the definition that relative Property (T) is
inherited by images. Consider the semidirect product $W\ltimes N$,
with group law
$(w_1,n_1)\cdot(w_2,n_2)=(w_1w_2,w_2^{-1}n_1w_2n_2)$. There is an
obvious morphism from $W\ltimes N$ to $W$ sending $(w,n)$ to $wn$.
Thus the lemma reduces to proving that $(W\ltimes N,N)$ has
relative Property (T).

Set $V=N\otimes_\Z\R$. Then $W\ltimes N$ naturally embeds as a
cocompact subgroup of cofinite volume into $W\ltimes V$. We are
going to show that $(W\ltimes V,V)$ has relative Property (T). It
trivially implies that $(W\ltimes V,N)$ has relative Property, and
then, since Property (T) relative to a given normal subgroup is
inherited by subgroups of cofinite volume \cite[Corollary
4.1(2)]{Jol}, this implies that $(W\ltimes V,N)$ has relative
Property (T).

So it remains to prove that $(W\ltimes V,V)$ has relative Property
(T). By a classical result by Burger, this reduces to proving that
the action by conjugation of $W\ltimes V$ on $V$ does not preserve
any Borel probability measure on the projective space $P(V^*)$, where $V$
denotes the dual space of $V$. This action factors through $W$,
so that we have to show that $W$ does not preserve any probability
on the projective space $P(V^*)$.

Let $V_1$ denote the vector subspace of $T_f$ generated by
$T_f\cap\alpha(W)$. We thus have to show that the action by
conjugation of $\alpha(W)$ on $P(T_f)$ does not preserve any
probability. Otherwise, by results of Furstenberg
\cite[§3.2]{Zimmer}, some finite index subgroup of $\alpha(W)$
preserve a nonzero subspace of $V_1^*$ on which it acts through
the action of a compact group. But this in contradication with the
fact that $\alpha(W)$ is Zariski dense in $\text{O}_f(B)$, and
that the action by conjugation of the connected semisimple group
without compact factors $\text{O}_0(\bar{B})$ on $T_f$ has no
invariant vectors.\epr

\bigskip
\footnotesize
\noindent Yves de Cornulier\\
E-mail: \url{yves.cornulier@math.u-psud.fr}



\end{document}